\documentclass[11pt,a4paper]{article}

\usepackage[usenames]{color}
\usepackage{amssymb}
\usepackage{amsmath}
\usepackage{amsthm}
\usepackage{amsfonts}
\usepackage{amscd}
\usepackage{graphicx}
\usepackage{mathtools}

\usepackage[colorlinks=true,
linkcolor=black,
filecolor=webbrown,
citecolor=black, urlcolor=blue]{hyperref}

\usepackage[cal=cm,scr=euler]{mathalpha}

\definecolor{webgreen}{rgb}{0,.5,0}
\definecolor{webbrown}{rgb}{.6,0,0}

\usepackage{color}
\usepackage{fullpage}
\usepackage{float}

\makeatletter
\newcommand{\defeq}{\mathrel{\aban@defeq}}
\newcommand{\aban@defeq}{%
  \vbox{\offinterlineskip\check@mathfonts
    \ialign{\hfil##\hfil\cr
      \fontsize{\ssf@size}{\z@}\normalfont def\cr
      \noalign{\kern1\p@}
      $\m@th=$\cr
      \noalign{\kern-.5\fontdimen22\textfont2}
    }%
  }%
}
\makeatother

\begin{document}

\def\refname {References}
\newtheorem{theorem} {Theorem}[section]
\newtheorem{lemma}[theorem] {Lemma}
\newtheorem{definition}[theorem] {Definition}
\newtheorem{example}[theorem] {Example}
\newtheorem{corollary}[theorem] {Corollary}
\newtheorem{proposition} [theorem]{Proposition}
\newtheorem{problem}{Problem}

\def\bibname {\bf References}

\title{\bf Existence and Uniqueness Property On a Generalized Ledin-Brousseau Sum}
\date{}
\author{Ivan Hadinata\footnote{Corresponding Author \\
2020 Mathematics Subject Classification: 11B39, 11B83.
\\
Keywords and phrases: Ledin-Brousseau sums, linear recurrence sequences, polynomial.
}}
\begingroup
\let\newpage\relax% Void the actions of \newpage
\maketitle
\endgroup
{\small \textbf{Abstract}. In this paper, we present the existence and uniqueness property on a finite sum involving a polynomial and a homogeneous linear recurrence sequence. This finite sum is of the form $\sum_{k=1}^n P(k)s_{hk+r}$ where $n$ is a positive integer, $P(x)$ is a polynomial in $\mathbb C[x]$, $h$ and $r$ are some integers, and $(s_k)_{k\in\mathbb Z}$ is a homogeneous linear recurrence sequence of degree $m\geq 2$ with some constraints. }

\section{Introduction} 

One of the frequently discussed topics involving finite sum and linear recurrence sequence is Ledin-Brousseau sum and its generalizations. In the beginning, Ledin \cite{ledin-1} and Brousseau \cite{brousseau-1} respectively studied the finite sums $\sum_{k=1}^{n}k^mF_k$ and $\sum_{k=1}^{n}k^mF_{k+r}$ where $(F_k)_{k\geq 0}$ is the classic Fibonacci sequence. Ledin \cite{ledin-1} presented the identity
\begin{equation}\label{S(m,n)}
    S(m,n) \defeq \sum_{k=1}^{n}k^mF_k = F_{n+1}P_2(m,n) + F_nP_1(m,n) + C(m)
\end{equation}
where $P_1(m,n)$ and $P_2(m,n)$ are polynomials in $n$ of degree $m$ and $C(m)$ is a constant depending on $m$. In the same year with Ledin, Brousseau \cite{brousseau-1} formulated the sum 
\begin{equation}\label{S_r(m,n)}
    S_r(m,n)\defeq \sum_{k=1}^{n-1}k^mF_{k+r} = \sum_{t=0}^{m}(-1)^t\Delta^t(n^m)F_{n+r+2t+1} + C
\end{equation}
where $C$ is a constant independent of $n$ and $\Delta$ is the forward difference operator with the properties $\Delta(f(n)) = f(n+1) - f(n)$, $\Delta^{t+1}(f(n)) = \Delta(\Delta^t(f(n)))$, and $\Delta^0(f(n)) = f(n)$.

Many various methods have been developed to evaluate the formula of $S(m,n)$ and $S_r(m,n)$. Ledin \cite{ledin-1} used an integration method to find the recursion of $P_1(m,n)$ and $P_2(m,n)$. Brousseau \cite{brousseau-1} used finite difference method to obtain the expression \eqref{S_r(m,n)}. Shannon and Ollerton, in their works \cite{ollerton-1} and \cite{shannon-1}, used the matrix method and formulated the linear recurrence relation of $S(m,n)$. Besides that, partial differentiation method was applied by Adegoke \cite{adegoke-1} to obtain the recursive relation for $P_1(m,n)$, $P_2(m,n), C(m)$ and then the polynomial form of $S(m,n)$. Lastly, Dresden \cite{dresden-1} evaluated $S(m,n)$ by applying a convolution representation of $F_n - n^m$ which involves binomial coefficients.

Ledin's and Brousseau's works (\cite{ledin-1, brousseau-1}) have also motivated many researchers afterward to investigate the more general sums. Mostly, the Ledin-Brousseau sum is generalized by replacing the Fibonacci sequence in the sums \eqref{S(m,n)} and \eqref{S_r(m,n)} with other linear recurrence sequences. Zeitlin \cite{zeitlin-1} extended the summations by using Gibonacci sequence. Adegoke \cite{adegoke-1} and Nair \cite{nair-1} generalized Ledin-Brousseau sum by replacing with Horadam sequence. Nair and Karunakaran \cite{nair-2} explored the Brousseau sum which is generated by k-Fibonacci numbers. Hadinata \cite{hadinata-1} considered the generalized Ledin-Brousseau sum by using a certain generalized second-degree linear recurrence sequence. Recently, the Brousseau sum for Tribonacci numbers was presented by Nair and Karunakaran \cite{nair-3}.

In this paper, we observe a generalized Ledin-Brousseau sum generated by the terms of a sequence defined in Definition \ref{definition of (s_k)}. Let, in all part of this paper, the notation $(s_k)_{k\in\mathbb Z}$ be referred to the sequence defined in Definition \ref{definition of (s_k)}.
\begin{definition}\label{definition of (s_k)}
Let $m\geq 2$ be a positive integer. Let $(s_k)_{k\in\mathbb Z}$ be a bi-infinite, homogeneous, and linear recurrence sequence of minimal degree $m$ which satisfy the following properties: 
\begin{itemize}
\item $(s_k)_{k\in\mathbb Z}$ satisfies the recurrence relation
\begin{equation}\label{recurrence of s_k}
    s_{k+m} = a_ms_{k+m-1} + a_{m-1}s_{k+m-2} + \cdots\cdots + a_2s_{k+1} + a_1s_k \quad (\forall k\in\mathbb Z)
\end{equation}
where $a_1, a_2, \ldots, a_m$ are complex numbers with $a_1\neq 0$. 
\item $(s_k)_{k\in\mathbb Z}$ does not satisfy all homogeneous linear recurrence relations of complex constant coefficients with degree $<m$. 
\item $s_1,s_2,\ldots,s_m$ are not simultaneously zero.
\item The characteristic polynomial of recurrence relation \eqref{recurrence of s_k}, that is,
\begin{equation}\label{characteristic polynomial}
    \mathcal{CP}(x) = x^m - a_mx^{m-1} - a_{m-1}x^{m-2} - \ldots - a_2x - a_1
\end{equation}
have $m$ distinct roots.
\end{itemize}
\end{definition}

\vspace{0.2cc}
We consider the more general Ledin-Brousseau sum $\sum_{k=1}^{n}P(k)s_{hk+r}$ where $n\in\mathbb N$, $P(x)\in\mathbb C[x]$ is a polynomial, and $h,r$ are integers with some constraints. We represent the sum $\sum_{k=1}^{n}P(k)s_{hk+r}$ in the form of $s_{(n+1)h+r},s_{(n+2)h+r},\ldots,s_{(n+m)h+r}$ and some polynomials. Then, we investigate the existence and uniqueness property behind this representation. There are two main results in this paper: Theorem \ref{theorem 1} and Theorem \ref{theorem 2}. Respectively, they present existence and uniqueness property about the representation of summation $\sum_{k=1}^n P(k)s_{hk+r}$.

\vspace{0.5cc}
\begin{theorem}\label{theorem 1}
    Let $P(x)\in\mathbb C[x]$ be a polynomial and $h,r$ be integers such that $h\neq 0$ and $h^{th}$ power of every roots of the characteristic polynomial \eqref{characteristic polynomial} are distinct and not equal to $1$. Then, there exist $m+1$ polynomials $P_1(x), P_2(x), \ldots, P_{m+1}(x)$ in $\mathbb C[x]$ such that the following identity holds:
    \begin{equation}\label{equation in theorem 1}
        \sum_{k=1}^{n}P(k)s_{hk+r} = \Big[\sum_{k=1}^{m}P_k(n)s_{(n+k)h+r} \Big]  + P_{m+1}(n), \;\forall n\in\mathbb N.
    \end{equation}
\end{theorem}

\vspace{1.5cc}
To show the existence of $(m+1)$-tuple of polynomials $(P_1(x), P_2(x), \ldots, P_{m+1}(x))\in\mathbb C[x]^{m+1}$ which satisfies \eqref{equation in theorem 1}, it is enough to find only one possible tuple satisfying \eqref{equation in theorem 1}. But, the harder problem to solve is "Is the tuple $(P_1(x), P_2(x), \ldots, P_{m+1}(x))\in\mathbb C[x]^{m+1}$ which satisfy \eqref{equation in theorem 1} unique or not ?". In this paper, we find that, for each $P(x)\in\mathbb C[x]$, the number of tuples $(P_1(x), P_2(x), \ldots, P_{m+1}(x))\in\mathbb C[x]^{m+1}$ which satisfy \eqref{equation in theorem 1} is always one.

\vspace{1cc}

\begin{theorem}\label{theorem 2}
    Let $P(x)\in\mathbb C[x]$ be a polynomial and $h,r$ be integers such that $h\neq 0$ and $h^{th}$ power of every roots of the characteristic polynomial \eqref{characteristic polynomial} are distinct and not equal to $1$. Then, the number of possible tuples $(P_1(x), P_2(x), \ldots, P_{m+1}(x))\in\mathbb C[x]^{m+1}$ which satisfy the identity \eqref{equation in theorem 1} is one. 
\end{theorem}

\vspace{1cc}

\section{Basic Results}\label{section 2}
\subsection{An Overview of $(s_k)_{k\in\mathbb Z}$}
Let $m$ roots of $\mathcal{CP}(x) = x^m - a_mx^{m-1} - a_{m-1}x^{m-2} - \ldots - a_2x - a_1$ be $r_1$, $r_2$, $\ldots$, $r_m$ $\in\mathbb C$ where all of $r_1,r_2,\ldots,r_m$ are distinct and $|r_1|\leq |r_2| \leq \cdots \leq |r_m|$. Since $a_1\neq 0$, then all of $r_1, r_2, \ldots, r_m$ are non-zero. By all properties in Definition \ref{definition of (s_k)}, the explicit formula of each $s_k$ is represented as 
\begin{equation}\label{explicit.formula.(s_k)}
    s_k = L_1r_1^k + L_2r_2^k + \cdots + L_mr_m^k, \quad\forall k\in\mathbb Z
\end{equation}
where $L_1$, $L_2$, $\ldots$, $L_m$ are some non-zero constants in $\mathbb C$. The values of $L_1$, $L_2$, $\ldots$, $L_m$ are obtained from system of $m$ equations $s_k = L_1r_1^k + L_2r_2^k + \cdots + L_mr_m^k$ $(k=1,2,\ldots,m)$. Hence, $(L_1, L_2,\ldots,L_m)$ is the solution of matrix equation
\begin{equation*}
    \begin{bmatrix}
        r_1 & r_2 & r_3 & \dots & r_m \\
        r_1^2 & r_2^2 & r_3^2 & \dots & r_m^2 \\
        r_1^3 & r_2^3 & r_3^3 & \dots & r_m^3 \\
        \vdots & \vdots & \vdots & \ddots & \vdots \\
        r_1^m & r_2^m & r_3^m & \dots & r_m^m
    \end{bmatrix}
    \begin{bmatrix}
        L_1 \\
        L_2 \\
        L_3 \\
        \vdots \\
        L_m
    \end{bmatrix} = 
    \begin{bmatrix}
        s_1 \\
        s_2 \\
        s_3 \\
        \vdots \\
        s_m
    \end{bmatrix}.
\end{equation*}
By using Cramer rule, for every $k\in\{1, 2, 3, \ldots, m\}$,
\begin{equation*}
    L_k = \frac{\det(A_k)}{\det\begin{bmatrix}
        r_1 & r_2 & r_3 & \dots & r_m \\
        r_1^2 & r_2^2 & r_3^2 & \dots & r_m^2 \\
        r_1^3 & r_2^3 & r_3^3 & \dots & r_m^3 \\
        \vdots & \vdots & \vdots & \ddots & \vdots \\
        r_1^m & r_2^m & r_3^m & \dots & r_m^m
    \end{bmatrix}} = \frac{\det(A_k)}{\prod_{i_1=1}^{m}r_{i_1} \times \prod_{1\leq i_1 < i_2 \leq m}(r_{i_2} - r_{i_1})},
\end{equation*}
where $A_k$ is the $m\times m$ matrix obtained by replacing the $k^{th}$ column of square matrix $\begin{bmatrix}
        r_1 & r_2 & r_3 & \dots & r_m \\
        r_1^2 & r_2^2 & r_3^2 & \dots & r_m^2 \\
        r_1^3 & r_2^3 & r_3^3 & \dots & r_m^3 \\
        \vdots & \vdots & \vdots & \ddots & \vdots \\
        r_1^m & r_2^m & r_3^m & \dots & r_m^m
    \end{bmatrix}$ with the column vector $\begin{bmatrix}
        s_1 \\
        s_2 \\
        s_3 \\
        \vdots \\
        s_m
    \end{bmatrix}$.

\vspace{0.5cc}
\subsection{Some Basic Lemmas}
In this subsection, we present some basic lemmas that will be used for showing Theorem \ref{theorem 1} and \ref{theorem 2}. 
\begin{lemma}\label{zero.polynom}
    Let $A$ be an infinite subset of $\mathbb N$. If $P(x)\in\mathbb C[x]$ is a polynomial and $\lim\limits_{\substack{n\to\infty \\ n\in A}}P(n)=0$, then $P(x)$ is a zero polynomial.
\end{lemma}
\begin{proof}
    If $P(x)$ is not a constant polynomial, then $P(x) = \delta_ux^u + \delta_{u-1}x^{u-1} + \cdots + \delta_1x+\delta_0$ for some $u\in\mathbb N$ and $\delta_0, \delta_1, \ldots, \delta_{u-1}, \delta_u \in \mathbb C$ with $\delta_u\neq 0$. We have
    \begin{equation*}
        0 = \lim\limits_{\substack{n\to\infty \\ n\in A}} |P(n)| = \lim\limits_{\substack{n\to\infty \\ n\in A}} |\delta_un^u + \delta_{u-1}n^{u-1} + \cdots + \delta_1n+\delta_0|,
    \end{equation*}
    but
    \begin{align*}
        \lim\limits_{\substack{n\to\infty \\ n\in A}} |\delta_un^u + \delta_{u-1}n^{u-1} + \cdots + \delta_1n+\delta_0| 
        &= \lim\limits_{\substack{n\to\infty \\ n\in A}} |\delta_un^u|\cdot\Big|1+\frac{\delta_{u-1}}{\delta_un} + \frac{\delta_{u-2}}{\delta_un^2} + \cdots + \frac{\delta_0}{\delta_un^u}\Big| \\
        &= \lim\limits_{\substack{n\to\infty \\ n\in A}} |\delta_un^u| \cdot \lim\limits_{\substack{n\to\infty \\ n\in A}} \Big|1+\frac{\delta_{u-1}}{\delta_un} + \frac{\delta_{u-2}}{\delta_un^2} + \cdots + \frac{\delta_0}{\delta_un^u}\Big| \\
        &= \infty\cdot 1 = \infty,
    \end{align*}
     so it is a contradiction and $P(x)$ can not be non-constant. \\
If $P(x)$ is a constant, let $P(x)\equiv c$, then the condition $\lim\limits_{\substack{n\to\infty \\ n\in A}}P(n)=0$ implies that $c=0$ and $P(x)\equiv 0$. \\
Hence $P(x)$ is a zero polynomial. 
\end{proof}

\vspace{0.5cc}
If $P(x)$ is a polynomial in $\mathbb C[x]$ and $y$ is a complex number with $|y|>1$, then the growth rate of $|y|^n$ is faster than $|P(n)|$ for all big enough positive integers $n$. It will cause that $P(n)/y^n \to 0$ as positive integer $n$ goes to $\infty$. This intuition motivates the following lemma:
\begin{lemma}\label{growth.rate.lemma}
    Let $A$ be an infinite subset of $\mathbb N$, $P(x)$ be a polynomial in $\mathbb C[x]$, and $y_1, y_2 \in \mathbb C$ be non-zero with $|y_1| > 1 > |y_2|$. Then,
    \begin{enumerate}
        \item $\lim\limits_{\substack{n\to\infty \\ n\in A}} \frac{P(n)}{y_1^n} = 0$. 
        \item $\lim\limits_{\substack{n\to\infty \\ n\in A}} P(n)y_2^n = 0$.
    \end{enumerate}
\end{lemma}
\begin{proof}
    It suffices to show the first limit. Because if the first limit is true, by setting $y_1:=1/y_2$ where $y_2\in\mathbb C$ is non-zero with $|y_2| < 1$, it yields the second limit. \\
    Here is proof of the first limit: Let $P(x) = \delta_ux^u + \delta_{u-1}x^{u-1} + \cdots + \delta_1x + \delta_0$ for some $u\in\mathbb N_0$ and $\delta_0, \delta_1, \ldots, \delta_u \in\mathbb C$. For all $n\in A$, we have
    \begin{align*}
        0 \leq \Big| \frac{P(n)}{y_1^n} \Big| = \frac{|\delta_un^u + \delta_{u-1}n^{u-1} + \cdots + \delta_1n + \delta_0|}{|y_1|^n} &\leq \frac{|\delta_u|n^u + |\delta_{u-1}|n^{u-1} + \cdots + |\delta_1|n + |\delta_0|}{|y_1|^n} \\
        &\leq \frac{M\cdot n^u}{|y_1|^n}
    \end{align*}
    where $M =|\delta_0| + |\delta_1| + \cdots + |\delta_u|$. \\
    Besides that, for all $n\in A$ with $n > \frac{(u+2)!}{(\ln{|y_1|})^{u+2}}$, we have
    \begin{equation*}
        0 < \frac{n^{u+1}}{|y_1|^n} = \frac{n^{u+1}}{\sum_{a=0}^{\infty}\frac{(n\ln{|y_1|})^a}{a!}} \leq \frac{n^{u+1}}{\frac{(n\ln|y_1|)^{u+2}}{(u+2)!}} = \frac{(u+2)!}{n(\ln |y_1|)^{u+2}} < 1  \;\implies\; 0 < \frac{n^u}{|y_1|^n} < \frac{1}{n}.
    \end{equation*}
    Since $\lim\limits_{\substack{n\to\infty \\ n\in A}} 0 = \lim\limits_{\substack{n\to\infty \\ n\in A}} \frac{1}{n} = 0$, by squeeze theorem, we get $\lim\limits_{\substack{n\to\infty \\ n\in A}} \frac{n^u}{|y_1|^n} = 0$. Then $\lim\limits_{\substack{n\to\infty \\ n\in A}} 0 = \lim\limits_{\substack{n\to\infty \\ n\in A}}\frac{M\cdot n^u}{|y_1|^n} = 0$. By squeeze theorem again, we obtain $\lim\limits_{\substack{n\to\infty \\ n\in A}} \Big| \frac{P(n)}{y_1^n} \Big| = 0$, then $\lim\limits_{\substack{n\to\infty \\ n\in A}} \frac{P(n)}{y_1^n} = 0$.
\end{proof}

\vspace{0.5cc}
The linear recurrence relation \eqref{recurrence of s_k} shows us that every term in $(s_k)_{k\in\mathbb Z}$ is representable as a fixed linear combination of $m$ previous terms. Then, the new question arises: "For a non-zero integer $\ell$, can $s_{m\ell + q}$ be represented as a linear combination of $s_{(m-1)\ell + q}, s_{(m-2)\ell + q}, \ldots , s_{\ell + q}, s_q$ for every integers $q$? If yes, are the coefficients of those linear combination fixed over all $q\in\mathbb Z$?". The answer is positive, by the constraint that $r_1^\ell, r_2^\ell, \ldots, r_m^\ell$ are distinct.
\begin{definition}
    For every non-zero integers $\ell$, let $e_1(\ell), e_2(\ell), \ldots, e_m(\ell)$ be complex numbers defined by
    \begin{equation*}
        e_{m+1-k}(\ell) = \sum_{1\leq i_1 < i_2 \cdots <i_k \leq m} (-1)^{k+1}(r_{i_1}r_{i_2}\cdots r_{i_k})^\ell, \;\;\forall k\in\{1,2,3,\ldots,m\}.
    \end{equation*}
\end{definition}
\vspace{0.5cc}
\begin{lemma}
    Let $\ell$ be a non-zero integer so that $r_1^\ell, r_2^\ell, \ldots, r_m^\ell$ are distinct, then for all $q\in\mathbb Z$:
    \begin{equation}\label{definition of e_i(.)}
        s_{m\ell+q} = e_m(\ell)s_{(m-1)\ell+q} + e_{m-1}(\ell)s_{(m-2)\ell+q} + \cdots + e_2(\ell)s_{\ell+q} + e_1(\ell)s_q.
    \end{equation}
\end{lemma}
\begin{proof}
    Let us consider an arbitrary $k'\in\{0,1,2,\ldots,\ell-1\}$. Let $T_t(k') = s_{k'+t\ell}$ for all $t\in\mathbb Z$. Then, for all $t\in\mathbb Z$,
    \begin{align*}
        T_t(k') = s_{k'+t\ell} &= L_1r_1^{k'+t\ell} + L_2r_2^{k'+t\ell} + \cdots + L_mrm^{k'+t\ell} \\
        &= L_1r_1^{k'}(r_1^\ell)^t + L_2r_2^{k'}(r_2^\ell)^t + \cdots + L_mr_m^{k'}(r_m^\ell)^t.
    \end{align*}
    Since $r_1^\ell, r_2^\ell,\cdots, r_m^\ell$ are distinct and $L_1r_1^{k'}, L_2r_2^{k'}, \cdots, L_mr_m^{k'}$ are non-zero constant, then the sequence $(T_t(k'))_{t\in\mathbb Z}$ satisfies a homogeneous linear recurrence relation of degree $m$ with constant coefficients. That is, there are $m$ complex constants $f_1, f_2, \cdots, f_m$ so that $f_1\neq 0$ and
    \begin{equation*}
        T_{t+m}(k') = f_mT_{t+m-1}(k') + f_{m-1}T_{t+m-2}(k') + \cdots + f_1T_t(k')
    \end{equation*}
    for all $t\in\mathbb Z$. This recurrence relation has the characteristic polynomial $x^m-f_mx^{m-1} - f_{m-1}x^{m-2} - \cdots - f_2x - f_1$ and the $m$ roots of this polynomial are $r_1^\ell, r_2^\ell, \cdots, r_m^\ell$. By Vieta's theorem, for all $k\in\{1,2,\cdots,m\}$,
    \begin{equation*}
        f_{m+1-k} = \sum_{1\leq i_1 < i_2 < \cdots < i_k \leq m} (-1)^{k+1}(r_{i_1}r_{i_2}\cdots r_{i_k})^\ell = e_{m+1-k}(\ell).
    \end{equation*}
    Therefore, for all $k'\in\{0,1,2,\ldots,\ell-1\}$ and all $t\in\mathbb Z$, 
    \begin{equation*}
        s_{k'+(t+m)\ell} = e_m(\ell)s_{k'+(t+m-1)\ell} + e_{m-1}(\ell)s_{k'+(t+m-2)\ell} + \cdots + e_2(\ell)s_{k'+(t+1)\ell} + e_1(\ell)s_{k'+t\ell}.
    \end{equation*}
    and it is equivalent to
    \begin{equation*}
        s_{m\ell+q} = e_m(\ell)s_{(m-1)\ell+q} + e_{m-1}(\ell)s_{(m-2)\ell+q} + \cdots + e_2(\ell)s_{\ell+q} + e_1(\ell)s_q \;\; (\forall q\in\mathbb Z)
    \end{equation*}
    as the lemma said.
\end{proof}

\vspace{0.5cc}
The next two lemmas present the properties involving the limit to infinity of a linear combination of exponential functions $\exp(...)$. 

\vspace{0.5cc}
\begin{lemma}\label{limit.sum.exp}
    Let $q$ be a positive integer, $c_1,c_2,\ldots,c_q$ be non-zero complex numbers, and $b_1,b_2,\ldots,b_q\in (0,2\pi)$ be distinct real numbers. Then, the limit
    \begin{equation*}
        \lim\limits_{\substack{n\to\infty \\ n\in\mathbb N}}(c_1\exp(inb_1)+c_2\exp(inb_2)+\cdots+c_q\exp(inb_q))
    \end{equation*}
    does not exist.
\end{lemma}
\begin{proof}
    Assume the contrary that the limit $\lim\limits_{\substack{n\to\infty \\ n\in\mathbb N}}(c_1\exp(inb_1)+c_2\exp(inb_2)+\cdots+c_q\exp(inb_q))$ exists, then there is a complex number $T$ such that $\lim\limits_{\substack{n\to\infty \\ n\in\mathbb N}}(c_1\exp(inb_1)+c_2\exp(inb_2)+\cdots+c_q\exp(inb_q)) = T$. It implies that for all $\alpha\in\mathbb Z$, we have
    \begin{equation}\label{useful.eq.12}
        \lim\limits_{\substack{n\to\infty \\ n\in\mathbb N}}(c_1\exp(i(n+\alpha)b_1)+c_2\exp(i(n+\alpha)b_2)+\cdots+c_q\exp(i(n+\alpha)b_q)) = T.
    \end{equation}
    By setting $\alpha=0,1,2,\ldots,q-1$ to \eqref{useful.eq.12}, we get the limit equation of matrix
    \begin{equation}\label{useful.eq.13}
        \lim\limits_{\substack{n\to\infty \\ n\in\mathbb N}}\mathcal{M}\begin{bmatrix}
            c_1e^{inb_1} \\
            c_2e^{inb_2} \\
            \vdots \\
            c_qe^{inb_q}
        \end{bmatrix} = 
        \begin{bmatrix}
            T \\
            T \\
            \vdots \\
            T
        \end{bmatrix}
    \end{equation}
    where $\mathcal{M}$ is $q\times q$ square matrix defined by
    \begin{equation*}
        \mathcal{M} = 
        \begin{bmatrix}
        1 & 1 & \dots & 1 \\
        e^{ib_1} & e^{ib_2} & \dots & e^{ib_q} \\
        (e^{ib_1})^2 & (e^{ib_2})^2 & \dots & (e^{ib_q})^2 \\
        \vdots & \vdots & \ddots & \vdots \\
        (e^{ib_1})^{q-1} & (e^{ib_2})^{q-1} & \dots & (e^{ib_q})^{q-1}
    \end{bmatrix}.
    \end{equation*}
    We observe that $\mathcal{M}$ is invertible since $\det(\mathcal{M}) = \prod_{1\leq i_1 < i_2 \leq q}(e^{ib_{i_2}}-e^{ib_{i_1}}) \neq 0$. Therefore, by equation \eqref{useful.eq.13}, we obtain
    \begin{equation*}
        \lim\limits_{\substack{n\to\infty \\ n\in\mathbb N}}\begin{bmatrix}
            c_1e^{inb_1} \\
            c_2e^{inb_2} \\
            \vdots \\
            c_qe^{inb_q}
        \end{bmatrix} = 
        \mathcal{M}^{-1}\begin{bmatrix}
            T \\
            T \\
            \vdots \\
            T
        \end{bmatrix} = 
        \begin{bmatrix}
            \mathscr{D}_1 \\
            \mathscr{D}_2 \\
            \vdots \\
            \mathscr{D}_q
        \end{bmatrix}
    \end{equation*}
    where $\mathscr{D}_1,\mathscr{D}_2,\ldots,\mathscr{D}_q$ are some complex numbers. Then, for all $k\in\{1,2,...,q\}$, $\lim\limits_{\substack{n\to\infty \\ n\in\mathbb N}}e^{inb_k} = \mathscr{D}_k/c_k$. It is an impossibility, because the condition $b_1,b_2,\ldots,b_q\in (0,2\pi)$ should imply that all the limits $\lim\limits_{\substack{n\to\infty \\ n\in\mathbb N}}e^{inb_1}$, $\lim\limits_{\substack{n\to\infty \\ n\in\mathbb N}}e^{inb_2}$, $\cdots$, $\lim\limits_{\substack{n\to\infty \\ n\in\mathbb N}}e^{inb_q}$ do not converge.
\end{proof}

\vspace{0.8cc}
 
\begin{lemma}\label{god.of.lemma}
    Let $q$ be a positive integer, $c_1,c_2,...,c_{q+1}$ be complex numbers, and $b_1,b_2,...,b_q\in (0,2\pi)$ be distinct real numbers. Then,
    \begin{equation*}
        \lim\limits_{\substack{n\to\infty \\ n\in\mathbb N}}(c_1\exp(inb_1)+c_2\exp(inb_2)+\cdots+c_q\exp(inb_q)) = c_{q+1}
    \end{equation*}
    if and only if $c_1=c_2=\cdots = c_{q+1} = 0$. 
\end{lemma}
\begin{proof}
    For the "if" part, it is trivial: If $c_1=c_2=\cdots=c_{q+1}=0$, then the limit equation $\lim\limits_{\substack{n\to\infty \\ n\in\mathbb N}}(c_1\exp(inb_1)+c_2\exp(inb_2)+\cdots+c_q\exp(inb_q)) = c_{q+1}$ holds.

    For the "only if" part: Let be given the limit equation $\lim\limits_{\substack{n\to\infty \\ n\in\mathbb N}}(c_1\exp(inb_1)+c_2\exp(inb_2)+\cdots+c_q\exp(inb_q)) = c_{q+1}$, then we must show $c_1=c_2=\cdots=c_{q+1}=0$. If one of $c_1,c_2,...,c_q$ is non-zero, then $\{k\in\{1,2,...,q\} : c_k \neq 0\} = \{\tau_1,\tau_2,...,\tau_j\}$ for some $j,\tau_1,\tau_2,...,\tau_j\in\{1,2,...,q\}$ with $\tau_1 < \tau_2 < \cdots < \tau_j$. It implies 
    \begin{equation}\label{useful.eq.14}
        \lim\limits_{\substack{n\to\infty \\ n\in\mathbb N}}(c_{\tau_1}\exp(inb_{\tau_1}) + c_{\tau_2}\exp(inb_{\tau_2}) + \cdots + c_{\tau_j}\exp(inb_{\tau_j})) = c_{q+1}.
    \end{equation}
    According to Lemma \ref{limit.sum.exp}, the limit $\lim\limits_{\substack{n\to\infty \\ n\in\mathbb N}}(c_{\tau_1}\exp(inb_{\tau_1})+c_{\tau_2}\exp(inb_{\tau_2})+\cdots+c_{\tau_j}\exp(inb_{\tau_j}))$ does not exist and this condition contradicts with \eqref{useful.eq.14}. Consequently, all of $c_1,c_2,...,c_q$ must be zero and it implies $c_{q+1}=0$.
\end{proof}

\vspace{1.5cc}

\section{Proof of Theorem \ref{theorem 1}}\label{section 3}
First, we define some notations used in this chapter as the following (definition \ref{defone-prooftheorem1} and \ref{deftwo-prooftheorem1}).
\begin{definition}\label{defone-prooftheorem1}
    For every $(d,\ell)\in\mathbb N_0 \times (\mathbb Z\setminus\{0\})$ so that $r_1^{\ell}, r_2^{\ell},\ldots,r_m^{\ell}$ are distinct and not equal to $1$, let $M(d,\ell) = (x_{i_1,i_2}(\ell))_{i_1,i_2=1}^{d+1}$ be a $(d+1)\times(d+1)$ matrix where $x_{i_1,i_2}(\ell)$ is the $i_1^{th}$-row and $i_2^{th}$-column entry of $M(d,\ell)$ defined by
    \begin{equation*}
        x_{i_1,i_2}(\ell) = \begin{dcases}
            \binom{i_2-1}{i_1-1}\sum_{i_3=1}^{m}e_i(\ell)(m+1-i_3)^{i_2-i_1}, \quad\textrm{if}\;\; i_1<i_2\\
            -1 + \sum_{i_3=1}^{m}e_{i_3}(\ell), \quad\textrm{if}\;\; i_1 = i_2\\
            0, \quad\textrm{if}\;\; i_1 > i_2.
        \end{dcases}
    \end{equation*}
    It can be checked that $M(d,\ell)$ is invertible with $det(M(d,\ell)) = (-1+\sum_{i_3=1}^{m}e_i(\ell))^{d+1} = (-\prod_{i_3=1}^{m}(1-r_{i_3}^\ell))^{d+1} \neq 0$. Then, the equation 
    \begin{equation*}
        M(d,\ell)\mathbf{v} = \left(\binom{d}{0}m^d \;\; \binom{d}{1}m^{d-1} \;\; \binom{d}{2}m^{d-2} \; \cdots \;\; \binom{d}{d}m^0 \right)^T
    \end{equation*}
    has a unique solution in $\mathbf{v}\in\mathbb C^{d+1}$. Let this unique solution be denoted by $(\delta_0(\ell) \;\; \delta_1(\ell) \;\; \cdots \;\; \delta_d(\ell))^T$.
\end{definition}

\vspace{0.5cc}
\begin{definition}\label{deftwo-prooftheorem1}
    For every $(d,\ell,r)\in\mathbb N_0\times (\mathbb Z\setminus\{0\}) \times \mathbb Z$ so that $r_1^\ell, r_2^\ell,\ldots, r_m^\ell$ are distinct and not equal to $1$, let the polynomials $P_{1,d,\ell,r}(x), P_{2,d,\ell,r}(x), \ldots , P_{m+1,d,\ell,r}(x)$ be defined as follows.
    \begin{equation*}
        P_{k,d,\ell,r}(x) = \begin{dcases}
            P_{m,d,\ell,r}(x-m+k) - \sum_{i_3=k+1}^{m}e_{i_3}(\ell)P_{m,d,\ell,r}(x+k+1-i_3), \quad\textrm{if}\;\; 1\leq k \leq m-1 \\
            \delta_d(\ell)x^d + \delta_{d-1}(\ell)x^{d-1} + \cdots + \delta_1(\ell)x + \delta_0(\ell), \quad\textrm{if}\;\; k=m \\
            s_{\ell +r} - \sum_{i_3=1}^{m}P_{i_3,d,\ell,r}(1)s_{(i_3+1)\ell +r}, \quad\textrm{if}\;\; k=m+1.
        \end{dcases}
    \end{equation*}
\end{definition}

\vspace{0.5cc}
\noindent
We then show the lemma for the special case of Theorem \ref{theorem 1} when $P(x)=x^d$ with $d\in\mathbb N_0$.
\begin{lemma}\label{lemma.in.prf.of.thm1.2}
    Given $(d,h,r)\in\mathbb N_0\times(\mathbb Z\setminus\{0\})\times\mathbb Z$ such that $r_1^h, r_2^h,\ldots, r_m^h$ are distinct and not equal to $1$, then
    \begin{equation*}
        \sum_{k=1}^{n} k^ds_{hk+r} = \Big[ \sum_{k=1}^{m}P_{k,d,h,r}(n)s_{(n+k)h+r} \Big] + P_{m+1,d,h,r}(n), \;\;\forall n\in\mathbb N.
    \end{equation*}
\end{lemma}
\begin{proof}
    We shall show the above identity by using induction. For all $n\in\mathbb N$, let $S(n)$ be the statement that $\sum_{k=1}^{n} k^ds_{hk+r} = \big[ \sum_{k=1}^{m}P_{k,d,h,r}(n)s_{(n+k)h+r} \big] + P_{m+1,d,h,r}(n)$. The definition of $P_{m+1,d,h,r}(x)$ implies that $S(1)$ is true. By considering an arbitrary $n_1\in\mathbb N$ which makes $S(n_1)$ true, we then have
    \begin{equation}\label{useful.eq.0}
        \sum_{k=1}^{n_1} k^ds_{hk+r} = \Big[ \sum_{k=1}^{m}P_{k,d,h,r}(n_1)s_{(n_1 +k)h+r} \Big] + P_{m+1,d,h,r}(n_1)
    \end{equation}
    By definition \ref{deftwo-prooftheorem1} and noting that $\sum_{i_3=m+1}^{m}e_{i_3}(h)P_{m,d,h,r}(x+m+1-i_3) = 0$, each $k\in\{2,3,\ldots,m\}$ satisfies
    \begin{equation*}
        P_{k,d,h,r}(x) = P_{m,d,h,r}(x-m+k) - \sum_{i_3=k+1}^{m}e_{i_3}(h)P_{m,d,h,r}(x+k+1-i_3)
    \end{equation*}
    and
    \begin{equation*}
        P_{k-1,d,h,r}(x+1) = P_{m,d,h,r}(x-m+k) - \sum_{i_3=k}^{m}e_{i_3}(h)P_{m,d,h,r}(x+k+1-i_3),
    \end{equation*}
    then
    \begin{equation}\label{useful.eq.1}
        P_{k,d,h,r}(x) - P_{k-1,d,h,r}(x+1) = e_k(h)P_{m,d,h,r}(x+1)
    \end{equation}
    Replacing $x$ in \eqref{useful.eq.1} with $x+m-k$ implies
    \begin{equation*}
        P_{k,d,h,r}(x+m-k) - P_{k-1,d,h,r}(x+m-k+1) = e_k(h)P_{m,d,h,r}(x+m-k+1).
    \end{equation*}
    Therefore, we obtain
    \begin{equation*}
        \sum_{k=2}^{m}(P_{k,d,h,r}(x+m-k) - P_{k-1,d,h,r}(x+m-k+1)) = \sum_{k=2}^{m} e_k(h)P_{m,d,h,r}(x+m-k+1)
    \end{equation*}
    \begin{equation}\label{useful.eq.2}
        \implies\;\; P_{m,d,h,r}(x) - P_{1,d,h,r}(x+m-1) = \sum_{k=2}^{m}e_k(h)P_{m,d,h,r}(x+m-k+1).
    \end{equation}
    Next, we observe that the given matrix equation
    \begin{equation*}
        M(d,h)\begin{bmatrix}
            \delta_0(h) \\
            \delta_1(h) \\
            \vdots \\
            \delta_d(h)
        \end{bmatrix} = \begin{bmatrix}
            \binom{d}{0}m^d \\
            \binom{d}{1}m^{d-1} \\
            \vdots \\
            \binom{d}{d}m^0
        \end{bmatrix}
    \end{equation*}
    is equivalent to the system of $d+1$ equations
    \begin{equation}\label{useful.eq.3}
        \sum_{i_1=1}^{m}\sum_{i_2=k}^{d}e_{i_1}(h)\delta_{i_2}(h)\binom{i_2}{k}(m+1-i_1)^{i_2-k} = \binom{d}{k}m^{d-k} + \delta_k(h), \;\; k=0,1,2,\ldots,d.
    \end{equation}
    Besides that, for all $k\in\{0,1,2,\ldots,d\}$, the coefficients of $x^k$ in the polynomials $(x+m)^d + P_{m,d,h,r}(x)$ and $\sum_{i_1=1}^{m}e_{i_1}(h)P_{m,d,h,r}(x+m+1-i_1)$ are respectively
    \begin{equation*}
        \binom{d}{k}m^{d-k} + \delta_k(h) \quad\;\textrm{and}\quad\; \sum_{i_1=1}^{m}\sum_{i_2=k}^{d}e_{i_1}(h)\delta_{i_2}(h)\binom{i_2}{k}(m+1-i_1)^{i_2-k}.
    \end{equation*}
    Since the degrees of $(x+m)^d + P_{m,d,h,r}(x)$ and $\sum_{i_3=1}^{m}e_{i_3}(h)P_{m,d,h,r}(x+m+1-i_3)$ are at most $d$, comparing to \eqref{useful.eq.3} yields
    \begin{equation}\label{useful.eq.4}
        (x+m)^d + P_{m,d,h,r}(x) = \sum_{i_3=1}^{m}e_{i_3}(h)P_{m,d,h,r}(x+m+1-i_3).
    \end{equation}
    We then decrease the equation \eqref{useful.eq.4} by \eqref{useful.eq.2} to obtain 
    \begin{equation}\label{useful.eq.5}
        (x+m)^d + P_{1,d,h,r}(x+m-1) = e_1(h)P_{m,d,h,r}(x+m).
    \end{equation}
    Replacing $x$ in \eqref{useful.eq.5} with $x-m+1$ yields the equation
    \begin{equation}\label{useful.eq.6}
        (x+1)^d + P_{1,d,h,r}(x) = e_1(h)P_{m,d,h,r}(x+1).
    \end{equation}
    Combining \eqref{useful.eq.1} and \eqref{useful.eq.6} implies that the ratio of polynomials $P_{m,d,h,r}(x+1) : (P_{m,d,h,r}(x)-P_{m-1,d,h,r}(x+1)) : (P_{m-1,d,h,r}(x)-P_{m-2,d,h,r}(x+1)) : \cdots\cdots : (P_{2,d,h,r}(x)-P_{1,d,h,r}(x+1)) : ((x+1)^d + P_{1,d,h,r}(x))$ is equal to $1 : e_m(h) : e_{m-1}(h) : \cdots : e_2(h) : e_1(h)$. By combining with the recurrence relation \eqref{definition of e_i(.)} when $\ell:=h$ and $q:=(n_1+1)h+r$, we will get the equation
    \begin{multline*}
        P_{m,d,h,r}(x+1)s_{(n_1+m+1)h+r} = ((x+1)^d+P_{1,d,h,r}(x))s_{(n_1+1)h+r} \\
        + \sum_{k=2}^{m}(P_{k,d,h,r}(x)-P_{k-1,d,h,r}(x+1))s_{(n_1+k)h+r}
    \end{multline*}
    which is equivalent to
    \begin{multline}\label{useful.eq.7}
        (x+1)^ds_{(n_1+1)h+r} = -P_{1,d,h,r}(x)s_{(n_1+1)h+r} + P_{m,d,h,r}(x+1)s_{n_1+m+1} \\
        + \sum_{k=2}^{m}(P_{k-1,d,h,r}(x+1)-P_{k,d,h,r}(x))s_{(n_1+k)h+r}.
    \end{multline}
    Setting $x:=n_1$ to \eqref{useful.eq.7} yields the equation
    \begin{multline}\label{useful.eq.8}
        (n_1+1)^ds_{(n_1+1)h+r} = -P_{1,d,h,r}(n_1)s_{(n_1+1)h+r} + P_{m,d,h,r}(n_1+1)s_{n_1+m+1} \\
        + \sum_{k=2}^{m}(P_{k-1,d,h,r}(n_1+1)-P_{k,d,h,r}(n_1))s_{(n_1+k)h+r}.
    \end{multline}
    By summing up the equations \eqref{useful.eq.0} and \eqref{useful.eq.8}, we obtain
    \begin{equation*}
        \sum_{k=1}^{n_1+1} k^ds_{hk+r} = \Big[ \sum_{k=1}^{m}P_{k,d,h,r}(n_1+1)s_{(n_1+k+1)h+r} \Big] + P_{m+1,d,h,r}(n_1+1).
    \end{equation*}
    Then, $S(n_1+1)$ is true. By induction, $S(n)$ is true for all $n\in\mathbb N$, then the proof (of Lemma \ref{lemma.in.prf.of.thm1.2}) completes.
\end{proof}

\vspace{0.8cc}
By considering a polynomial $P(x)\in\mathbb C[x]$, we have $P(x) =\sum_{d\in\mathbb N_0} Coef(d,P(x))x^d$ where $Coef(d,P(x))$ is the coefficient of $x^d$ in $P(x)$. Then, a suitable example of polynomials $P_1(x), P_2(x),\ldots,P_{m+1}(x)$ in $\mathbb C[x]$ which satisfy the identity \eqref{equation in theorem 1} is
\begin{equation*}
    P_k(x) := \mathscr{P}_k(x) \defeq \sum_{d\in\mathbb N_0}Coef(d,P(x))P_{k,d,h,r}(x), \;\;\forall k\in\{1,2,\ldots,m+1\}.
\end{equation*}
Hence, the proof (of Theorem \ref{theorem 1}) completes.

\vspace{2cc}
\section{Proof of Theorem \ref{theorem 2}}\label{section 4}
\noindent
Let us consider a polynomial $P(x)$ in $\mathbb C[x]$. Next, we suppose that $\mathscr{S}$ is the set of $(m+1)$-tuple of polynomials $(\gamma_1(x),\gamma_2(x),\ldots,\gamma_{m+1}(x))\in\mathbb C[x]^{m+1}$ which satisfy
    \begin{equation}\label{eq.firstlemmasect4}
        \gamma_1(n)s_{(n+1)h+r} + \gamma_2(n)s_{(n+2)h+r} + \cdots + \gamma_m(n)s_{(n+m)h+r} + \gamma_{m+1}(n) = 0 \quad (\forall n\in\mathbb N)
    \end{equation}
    and $\mathscr{S}(P(x))$ is the set of all possible $(P_1(x), P_2(x),\ldots,P_{m+1}(x))$ in $\mathbb C[x]^{m+1}$ which satisfy the identity \eqref{equation in theorem 1}.

Then, we have an important lemma as follows.
\begin{lemma}\label{firstlemmasect4}
     $|\mathscr{S}(P(x))|=|\mathscr{S}|$.
\end{lemma}
\begin{proof}
    As we have said in the last paragraph of Section \ref{section 3}, $(\mathscr{P}_1(x), \mathscr{P}_2(x),\ldots,\mathscr{P}_{m+1}(x))\in\mathscr{S}(P(x))$. Let us observe that for every $(P_1(x), P_2(x),\ldots,P_{m+1}(x))$ in $\mathscr{S}(P(x))$, then $(P_1(x)-\mathscr{P}_1(x),P_2(x)-\mathscr{P}_2(x),\ldots,P_{m+1}(x)-\mathscr{P}_{m+1}(x)) \in \mathscr{S}$. So, we can construct a mapping $T:\mathscr{S}(P(x)) \to \mathscr{S}$ defined by
    \begin{equation*}
        T(P_1(x), P_2(x),\ldots,P_{m+1}(x)) = (P_1(x)-\mathscr{P}_1(x),P_2(x)-\mathscr{P}_2(x),\ldots,P_{m+1}(x)-\mathscr{P}_{m+1}(x))
    \end{equation*}
    for all $(P_1(x), P_2(x),\ldots,P_{m+1}(x))\in\mathscr{S}(P(x))$. \\
    $T$ is injective because if 
    \begin{equation*}
        \textrm{$(P_1^{(1)}(x),P_2^{(1)}(x),\cdots,P_{m+1}^{(1)}(x))\in\mathscr{S}(P(x))$ and $(P_1^{(2)}(x),P_2^{(2)}(x),\cdots,P_{m+1}^{(1)}(x))\in \mathscr{S}(P(x))$}
    \end{equation*}
    satisfy the condition 
    \begin{equation*}
        T(P_1^{(1)}(x),P_2^{(1)}(x),\cdots,P_{m+1}^{(1)}(x)) = T(P_1^{(2)}(x),P_2^{(2)}(x),\cdots,P_{m+1}^{(1)}(x)),
    \end{equation*}
    then for all $k\in\{1,2,...,m+1\}$,
    \begin{equation*}
        P_k^{(1)}(x) - \mathscr{P}_k(x) = P_k^{(2)}(x) - \mathscr{P}_k(x)  \;\;\implies\;\; P_k^{(1)}(x) = P_k^{(2)}(x),
    \end{equation*}
    hence $(P_1^{(1)}(x),P_2^{(1)}(x),\cdots,P_{m+1}^{(1)}(x)) = (P_1^{(2)}(x),P_2^{(2)}(x),\cdots,P_{m+1}^{(2)}(x))$. \\
    Besides that, $T$ is surjective because for every $(P_1^{(3)}(x),P_2^{(3)}(x),\cdots,P_{m+1}^{(3)}(x))\in \mathscr{S}$, there exists $(P_1^{(3)}(x)+\mathscr{P}_1(x),P_2^{(3)}(x)+\mathscr{P}_2(x),\cdots,P_{m+1}^{(3)}(x)+\mathscr{P}_{m+1}(x)) \in \mathscr{S}(P(x))$ such that
    \begin{equation*}
        T(P_1^{(3)}(x)+\mathscr{P}_1(x),P_2^{(3)}(x)+\mathscr{P}_2(x),\cdots,P_{m+1}^{(3)}(x)+\mathscr{P}_{m+1}(x)) = (P_1^{(3)}(x),P_2^{(3)}(x),\cdots,P_{m+1}^{(3)}(x)).
    \end{equation*}
    Thus, $T$ is a bijection and the result follows.
\end{proof}

By replacing each of $s_{(n+1)h+r}, s_{(n+2)h+r},\ldots, s_{(n+m)h+r}$ in \eqref{eq.firstlemmasect4} with the formula \eqref{explicit.formula.(s_k)}, the identity \eqref{eq.firstlemmasect4} is equivalent to
\begin{equation}\label{useful.eq.9}
    \gamma_{m+1}(n) + \sum_{i_1=1}^{m}\sum_{i_2=1}^{m}\gamma_{i_1}(n)L_{i_2}r_{i_2}^{(n+i_1)h+r} = 0 \;\quad (\forall n\in\mathbb N).
\end{equation}
By defining the polynomials $Q_1(x), Q_2(x), \ldots, Q_m(x)$ with 
\begin{equation}\label{useful.eq.33}
    Q_{i_2}(x) = L_{i_2}r_{i_2}^{h+r}\sum_{i_1=1}^{m}\gamma_{i_1}(x)r_{i_2}^{(i_1-1)h}, \quad\forall i_2=1,2,...,m,
\end{equation}
the identity \eqref{useful.eq.9} is equivalent to
\begin{equation}\label{useful.eq.10}
    Q_1(n)r_1^{hn} + Q_2(n)r_2^{hn} + \cdots + Q_m(n)r_m^{hn} + \gamma_{m+1}(n) = 0, \;\forall n\in\mathbb N.
\end{equation}

We shall show that all of $Q_1(x), Q_2(x), \ldots, Q_m(x)$ are zero polynomials. Assume the contrary that there exists a polynomial among $Q_1(x), Q_2(x), \ldots, Q_m(x)$ which is not a zero polynomial. Let $\{t_1,t_2,\cdots,t_p\}$, with $p\in\{1,2,...,m\}$ and $t_1<t_2<\cdots<t_p$, be the set of all indices $i_2$ in $\{1,2,...,m\}$ for which $Q_{i_2}(x)$ is not zero polynomial. Then, the identity \eqref{useful.eq.10} becomes
\begin{equation}\label{useful.eq.11}
    \gamma_{m+1}(n) + \sum_{k=1}^{p}Q_{t_k}(n)r_{t_k}^{hn} = 0,  \;\forall n\in\mathbb N.
\end{equation}
Suppose that $\{|r_{t_k}| : k\in\{1,2,...,p\} \} = \{\rho_1,\rho_2,...,\rho_w\}$ for some $w\in\mathbb N$ with $w\leq p$ and $0< \rho_1 < \rho_2 < \cdots < \rho_w$. Therefore, we can construct a function $\kappa:\{0,1,2,...,w\}\to\{0,1,2,...,p\}$ in such a way that $0=\kappa(0)<\kappa(1)<\cdots < \kappa(w) = p$ and for every $i_1\in\{1,2,...,w\}$,
\begin{equation*}
    \rho_{i_1} = \big|r_{t_{\kappa(i_1-1)+i_2}}^h \big|, \quad\forall i_2\in\{1,2,...,\kappa(i_1)-\kappa(i_1-1)\}.
\end{equation*}
Without loss of generality, for every $i_1\in\{1,2,...,w\}$ and $i_2\in\{1,2,...,\kappa(i_1)-\kappa(i_1-1)\}$, let $r_{t_{\kappa(i_1-1)+i_2}} = \rho_{i_1}\exp(i\theta_{i_1,i_2})$ where $0\leq \theta_{i_1,1}<\theta_{i_1,2}<\cdots<\theta_{i_1,\kappa(i_1)-\kappa(i_1-1)} < 2\pi$ ,and let the polynomials $Q_{i_1,i_2}(x) = Q_{t_{\kappa(i_1-1)+i_2}}(x)$. Then, the identity \eqref{useful.eq.11} is equivalent to
\begin{equation}\label{useful.eq.20}
    \gamma_{m+1}(n) + \sum_{i_1=1}^{w}\rho_{i_1}^nV_n(i_1) = 0, \;\;\forall n\in\mathbb N,
\end{equation}
where $V_n(i_1) = \sum_{i_2=1}^{\kappa(i_1)-\kappa(i_1-1)}Q_{i_1,i_2}(n)\exp(in\theta_{i_1,i_2})$ for every $(i_1,n)\in\{1,2,...,w\}\times\mathbb N$. 

\vspace{1cc}
\noindent
Let us consider the identity \eqref{useful.eq.20} into 3 cases: $\rho_w <1$, $\rho_w > 1$, and $\rho_w = 1$.

\vspace{1cc}
\noindent
\underline{CASE $\rho_w < 1$}: \\
We note that $\rho_1,\rho_2,...,\rho_w <1$. By Lemma \ref{growth.rate.lemma}, for each $i_1\in\{1,2,...,w\}$ and $i_2\in\{1,2,...,\kappa(i_1)-\kappa(i_1-1)\}$, we have
\begin{equation*}
    \lim\limits_{\substack{n\to\infty \\ n\in\mathbb N}} |\rho_{i_1}^nQ_{i_1,i_2}(n)\exp(in\theta_{i_1,i_2})| = \lim\limits_{\substack{n\to\infty \\ n\in\mathbb N}} |\rho_{i_1}^nQ_{i_1,i_2}(n)| = 0 \;\implies\; \lim\limits_{\substack{n\to\infty \\ n\in\mathbb N}} \rho_{i_1}^nQ_{i_1,i_2}(n)\exp(in\theta_{i_1,i_2}) = 0. 
\end{equation*}
Limiting $n\to\infty$ to \eqref{useful.eq.20} yields $\lim\limits_{\substack{n\to\infty \\ n\in\mathbb N}}\gamma_{m+1}(n) = 0$. By Lemma \ref{zero.polynom}, $\gamma_{m+1}(x)$ is a zero polynomial. Then, the identity \eqref{useful.eq.20} is equivalent to
\begin{equation}\label{useful.eq.25}
    \sum_{i_1=1}^{w}\frac{\rho_{i_1}^n}{\rho_w^n}V_n(i_1) = 0, \quad\forall n\in\mathbb N.
\end{equation}
If $w=1$, then the identity \eqref{useful.eq.25} yields $V_n(w)=0$ for all $n\in\mathbb N$, so we can write $\lim\limits_{\substack{n\to\infty \\ n\in\mathbb N}} V_n(w) = 0$. If $w>1$, then for all $i_1\in\{1,2,...,w-1\}$, $\frac{\rho_{i_1}}{\rho_w} < 1$. By triangle inequality and Lemma \ref{growth.rate.lemma}, each $i_1$ in $\{1,2,...,w-1\}$ satisfies 
\begin{equation*}
    \lim\limits_{\substack{n\to\infty \\ n\in\mathbb N}}\bigg|\frac{\rho_{i_1}^n}{\rho_w^n}V_n(i_1)\bigg| \leq \lim\limits_{\substack{n\to\infty \\ n\in\mathbb N}}\sum_{i_2 = 1}^{\kappa(i_1)-\kappa(i_1-1)}\bigg|\frac{\rho_{i_1}^n}{\rho_w^n}Q_{i_1,i_2}(n)\bigg| = 0,
\end{equation*}
hence $\lim\limits_{\substack{n\to\infty \\ n\in\mathbb N}} \frac{\rho_{i_1}^n}{\rho_w^n}V_n(i_1) = 0$. Limiting $n\to\infty$ to \eqref{useful.eq.25} implies $\lim\limits_{\substack{n\to\infty \\ n\in\mathbb N}} V_n(w)=0$. \\
Now, we focus on the limit $\lim\limits_{\substack{n\to\infty \\ n\in\mathbb N}} V_n(w) = 0$. It is clear that there is a $w'\in\{1,2,...,\kappa(w)-\kappa(w-1)\}$ in such a way that $Q_{w,w'}(x)$ has the highest degree among the polynomials $Q_{w,1}(x),Q_{w,2}(x),...,Q_{w,w'}(x)$. Then, $\lim\limits_{\substack{n\to\infty \\ n\in\mathbb N}} \frac{1}{Q_{w,w'}(n)}$ equals zero or a complex constant. It implies the following limit.
\begin{equation}\label{useful.eq.26}
     \lim\limits_{\substack{n\to\infty \\ n\in\mathbb N}} \sum_{i_2=1}^{\kappa(w)-\kappa(w-1)}\frac{Q_{w,i_2}(n)}{Q_{w,w'}(n)}\exp(in\theta_{w,i_2}) = \lim\limits_{\substack{n\to\infty \\ n\in\mathbb N}} \frac{V_n(w)}{Q_{w,w'}(n)} = 0
\end{equation}
For all $i_2\in\{1,2,...,\kappa(w)-\kappa(w-1)\}$, let us define $\mathcal{J}_{i_2} \defeq \lim\limits_{\substack{n\to\infty \\ n\in\mathbb N}} \frac{Q_{w,i_2}(n)}{Q_{w,w'}(n)}$. Then, each $i_2$ in $\{1,2,...,\kappa(w)-\kappa(w-1)\}$ satisfies the limit
\begin{equation*}
    \lim\limits_{\substack{n\to\infty \\ n\in\mathbb N}}\left(\mathcal{J}_{i_2}-\frac{Q_{w,i_2}(n)}{Q_{w,w'}(n)} \right)\exp(in\theta_{w,i_2}) = 0,
\end{equation*}  
which derives from
\begin{equation*}
    \lim\limits_{\substack{n\to\infty \\ n\in\mathbb N}}\bigg| \left(\mathcal{J}_{i_2}-\frac{Q_{w,i_2}(n)}{Q_{w,w'}(n)} \right)\exp(in\theta_{w,i_2}) \bigg| = \lim\limits_{\substack{n\to\infty \\ n\in\mathbb N}}\bigg| \mathcal{J}_{i_2}-\frac{Q_{w,i_2}(n)}{Q_{w,w'}(n)}  \bigg| = |\mathcal{J}_{i_2}-\mathcal{J}_{i_2}| = 0.
\end{equation*}
Consequently,
\begin{equation}\label{useful.eq.27}
    \lim\limits_{\substack{n\to\infty \\ n\in\mathbb N}} \sum_{i_2=1}^{\kappa(w)-\kappa(w-1)}\left(\mathcal{J}_{i_2}-\frac{Q_{w,i_2}(n)}{Q_{w,w'}(n)} \right)\exp(in\theta_{w,i_2}) = 0.
\end{equation}
Summing up the limit equation \eqref{useful.eq.26} and \eqref{useful.eq.27} yields
\begin{equation}\label{useful.eq.28}
    \lim\limits_{\substack{n\to\infty \\ n\in\mathbb N}}\left( \mathcal{J}_1\exp(in\theta_{w,1}) + \cdots + \mathcal{J}_{\kappa(w)-\kappa(w-1)}\exp(in\theta_{w,\kappa(w)-\kappa(w-1)}) \right) = 0
\end{equation}
If $\kappa(w)-\kappa(w-1) = 1$, we have $w'=1$ and 
\begin{equation*}
    \lim\limits_{\substack{n\to\infty \\ n\in\mathbb N}}|\mathcal{J}_1| = \lim\limits_{\substack{n\to\infty \\ n\in\mathbb N}} |\mathcal{J_1}\exp(in\theta_{w,1})| = 0 \;\implies\; \mathcal{J}_1 = 0,
\end{equation*}
which contradicts to the fact $\mathcal{J}_{w'} = \mathcal{J}_1 = 1$. \\
If $\kappa(w)-\kappa(w-1) > 1$, by \eqref{useful.eq.28} and using Lemma \ref{god.of.lemma}, we must have $\mathcal{J}_1 = \mathcal{J}_2 = \cdots = \mathcal{J}_{\kappa(w)-\kappa(w-1)} = 0$. It also contradicts to the fact $\mathcal{J}_{w'} = 1$.

\vspace{1cc}
\noindent
\underline{CASE $\rho_w > 1$}: \\
The identity \eqref{useful.eq.20} is equivalent to
\begin{equation}\label{useful.eq.29}
    \frac{\gamma_{m+1}(n)}{\rho_w^n} + \sum_{i_1=1}^{w}\frac{\rho_{i_1}^n}{\rho_w^n}V_n(i_1) = 0, \quad\forall n\in\mathbb N.
\end{equation}
By Lemma \ref{growth.rate.lemma}, $\lim\limits_{\substack{n\to\infty \\ n\in\mathbb N}}\frac{\gamma_{m+1}(n)}{\rho_w^n} = 0$. \\
If $w=1$, by limiting $n\to\infty$ to \eqref{useful.eq.29}, we get $\lim\limits_{\substack{n\to\infty \\ n\in\mathbb N}} V_n(w) = 0$. \\
If $w>1$, then $\frac{\rho_{i_2}}{\rho_w} < 1$ for all $i_2\in\{1,2,...,w-1\}$. By triangle inequality and Lemma \ref{growth.rate.lemma}, each $i_1$ in $\{1,2,...,w-1\}$ satisfies 
\begin{equation*}
    \lim\limits_{\substack{n\to\infty \\ n\in\mathbb N}} \bigg|\frac{\rho_{i_1}^n}{\rho_w^n}V_n(i_1)\bigg| \leq \lim\limits_{\substack{n\to\infty \\ n\in\mathbb N}}\sum_{i_2 = 1}^{\kappa(i_1)-\kappa(i_1-1)}\bigg|\frac{\rho_{i_1}^n}{\rho_w^n}Q_{i_1,i_2}(n)\bigg| = 0,
\end{equation*}
then $\lim\limits_{\substack{n\to\infty \\ n\in\mathbb N}}\frac{\rho_{i_1}^n}{\rho_w^n}V_n(i_1) = 0$. By limiting $n\to\infty$ to \eqref{useful.eq.29}, we obtain $\lim\limits_{\substack{n\to\infty \\ n\in\mathbb N}}V_n(w) = 0$. \\
The implication of the limit $\lim\limits_{\substack{n\to\infty \\ n\in\mathbb N}}V_n(w) = 0$ is as similar as our step in the case when $\rho_w < 1$. There will be a contradiction.

\vspace{1cc}
\noindent
\underline{CASE $\rho_w = 1$}: \\
If $w=1$, the identity \eqref{useful.eq.20} is equivalent to $V_n(w) + \gamma_{m+1}(n) = 0$, $\forall n\in\mathbb N$, hence we get the limit $\lim\limits_{\substack{n\to\infty \\ n\in\mathbb N}}(V_n(w)+\gamma_{m+1}(n)) = 0$. \\
If $w>1$, we have $\rho_{i_1} < 1$ for all $i_1\in\{1,2,...,w-1\}$. By triangle inequality and Lemma \ref{growth.rate.lemma}, for all $i_1\in\{1,2,...,w-1\}$, 
\begin{equation*}
    \lim\limits_{\substack{n\to\infty \\ n\in\mathbb N}}\big|\rho_{i_1}^nV_n(i_1)\big| \leq \lim\limits_{\substack{n\to\infty \\ n\in\mathbb N}}\sum_{i_2=1}^{\kappa(i_1)-\kappa(i_1-1)}\big|\rho_{i_1}^n Q_{i_1,i_2}(n)\big| = 0,
\end{equation*}
hence $\lim\limits_{\substack{n\to\infty \\ n\in\mathbb N}}\rho_{i_1}^nV_n(i_1) = 0$. Therefore, by limiting $n\to\infty$ to \eqref{useful.eq.20}, we get $\lim\limits_{\substack{n\to\infty \\ n\in\mathbb N}}(V_n(w)+\gamma_{m+1}(n)) = 0$. \\
Suppose $\mathcal{Q}(x)$ is the polynomial among $Q_{w,1}(x), Q_{w,2}(x), \ldots, Q_{w,\kappa(w)-\kappa(w-1)}(x), \gamma_{m+1}(x)$ which has the highest degree, so $\lim\limits_{\substack{n\to\infty \\ n\in\mathbb N}}\frac{1}{\mathcal{Q}(n)}$ equals zero or a complex constant. We obtain
\begin{equation}\label{useful.eq.30}
\lim\limits_{\substack{n\to\infty \\ n\in\mathbb N}}\left(\frac{\gamma_{m+1}(n)}{\mathcal{Q}(n)} + \sum_{i_2=1}^{\kappa(w)-\kappa(w-1)}\frac{Q_{w,i_2}(n)}{\mathcal{Q}(n)}\exp(in\theta_{w,i_2})\right) = \lim\limits_{\substack{n\to\infty \\ n\in\mathbb N}}\frac{V_n(w)+\gamma_{m+1}(n)}{\mathcal{Q}(n)} = 0. 
\end{equation}
Next, suppose $\; \mathscr{J}_0 \defeq\lim\limits_{\substack{n\to\infty \\ n\in\mathbb N}}\frac{\gamma_{m+1}(n)}{\mathcal{Q}(n)} \;$ and $\; \mathscr{J}_{i_2}\defeq\lim\limits_{\substack{n\to\infty \\ n\in\mathbb N}}\frac{Q_{w,i_2}(n)}{\mathcal{Q}(n)} \;$ for every $i_2\in\{1,2,...,\kappa(w)-\kappa(w-1)\}$. Let us observe that, for all $i_2\in\{1,2,...,\kappa(w)-\kappa(w-1)\}$,
\begin{equation*}
    \lim\limits_{\substack{n\to\infty \\ n\in\mathbb N}}\bigg|\left(\mathscr{J}_{i_2}-\frac{Q_{w,i_2}(n)}{\mathcal{Q}(n)}\right)\exp(in\theta_{w,i_2})\bigg| = \lim\limits_{\substack{n\to\infty \\ n\in\mathbb N}}\bigg|\mathscr{J}_{i_2}-\frac{Q_{w,i_2}(n)}{\mathcal{Q}(n)}\bigg| = \big|\mathscr{J}_{i_2}-\mathscr{J}_{i_2} \big| = 0
\end{equation*}
and therefore,
\begin{equation*}
    \lim\limits_{\substack{n\to\infty \\ n\in\mathbb N}}\left(\mathscr{J}_{i_2}-\frac{Q_{w,i_2}(n)}{\mathcal{Q}(n)}\right)\exp(in\theta_{w,i_2}) = 0.
\end{equation*}
So, we get the following limit of summation:
\begin{equation}\label{useful.eq.31}
    \lim\limits_{\substack{n\to\infty \\ n\in\mathbb N}}\sum_{i_2=1}^{\kappa(w)-\kappa(w-1)}\left(\mathscr{J}_{i_2}-\frac{Q_{w,i_2}(n)}{\mathcal{Q}(n)}\right)\exp(in\theta_{w,i_2}) = 0.
\end{equation}
Summing up the limit equations \eqref{useful.eq.30} and \eqref{useful.eq.31} implies the new limit equation
\begin{equation}\label{useful.eq.32}
    \lim\limits_{\substack{n\to\infty \\ n\in\mathbb N}}\left(\mathscr{J}_0 + \sum_{i_2=1}^{\kappa(w)-\kappa(w-1)}\mathscr{J}_{i_2}\exp(in\theta_{w,i_2})\right) = 0.
\end{equation}
Since $r_1^h,r_2^h,...,r_m^h$ are not equal to $1$ and $\rho_w = 1$, then we must have $0<\theta_{w,1}<\theta_{w,2}<\cdots < \theta_{w,\kappa(w)-\kappa(w-1)}<2\pi$. By applying Lemma \ref{god.of.lemma} to \eqref{useful.eq.32}, we obtain $\mathscr{J}_0 = \mathscr{J}_1 = \cdots = \mathscr{J}_{\kappa(w)-\kappa(w-1)} = 0$. It is clearly a contradiction because one of $\mathscr{J}_0,\mathscr{J}_1,...,\mathscr{J}_{\kappa(w)-\kappa(w-1)}$ should be equal to $1$.

Finally, all of $Q_1(x),Q_2(x),...,Q_m(x)$ are zero polynomials. By \eqref{useful.eq.10} and Lemma \ref{zero.polynom}, $\gamma_{m+1}(x)$ is a zero polynomial. By \eqref{useful.eq.33}, we obtain the following matrix equation.
\begin{equation}\label{useful.eq.34}
    \begin{bmatrix}
        1 & r_1^h & r_1^{2h} & \dots & r_1^{(m-1)h} \\
        1 & r_2^h & r_2^{2h} & \dots & r_2^{(m-1)h} \\
        1 & r_3^h & r_3^{2h} &\dots & r_3^{(m-1)h} \\
        \vdots & \vdots & \vdots & \ddots & \vdots \\
        1 & r_m^h & r_m^{2h} & \dots & r_m^{(m-1)h}
    \end{bmatrix}
    \begin{bmatrix}
        \gamma_1(x) \\
        \gamma_2(x) \\
        \gamma_3(x) \\
        \vdots \\
        \gamma_m(x)
    \end{bmatrix} = 
    \begin{bmatrix}
        0 \\
        0 \\
        0 \\
        \vdots \\
        0
    \end{bmatrix}.
\end{equation}
The matrix 
\begin{equation*}
    \begin{bmatrix}
        1 & r_1^h & r_1^{2h} & \dots & r_1^{(m-1)h} \\
        1 & r_2^h & r_2^{2h} & \dots & r_2^{(m-1)h} \\
        1 & r_3^h & r_3^{2h} &\dots & r_3^{(m-1)h} \\
        \vdots & \vdots & \vdots & \ddots & \vdots \\
        1 & r_m^h & r_m^{2h} & \dots & r_m^{(m-1)h}
    \end{bmatrix}
\end{equation*}
is invertible, because its determinant is $\prod_{1\leq i_1 < i_2 \leq m}(r_{i_2}^h - r_{i_1}^h)\neq 0$. Therefore, the equation \eqref{useful.eq.34} implies that all of $\gamma_1(x)$, $\gamma_2(x)$, $\ldots$, $\gamma_m(x)$ are zero polynomials. In conclusion, $\mathscr{S}$ has only one element. By Lemma \ref{firstlemmasect4}, $\mathscr{S}(P(x))$ has one element too. It means that there is uniquely $(m+1)$-tuple $(\gamma_1(x),\gamma_2(x),\ldots,\gamma_{m+1}(x))\in\mathbb C[x]^{m+1}$ which satisfies the identity
\begin{equation*}
    \sum_{k=1}^{n}P(k)s_{hk+r} = \Big[\sum_{k=1}^{m}P_k(n)s_{(n+k)h+r} \Big]  + P_{m+1}(n), \;\forall n\in\mathbb N,     
\end{equation*}
and the proof completes.

\vspace{1.5cc}
\noindent
\textbf{Ivan Hadinata}\\
Faculty of Mathematics and Natural Sciences, 
Gadjah Mada University, Yogyakarta, Indonesia\\
\textit{Email}: \href{mailto: ivanhadinata2005@mail.ugm.ac.id}{\tt ivanhadinata2005@mail.ugm.ac.id}

\end{document}